\title{ ~~\\ On Robin's criterion for the\\ Riemann Hypothesis}
\author{Y.-J.~Choie, N.~Lichiardopol, P.~Moree, P.~Sol\'e}
\documentclass[12pt]{article}
\usepackage{amssymb, latexsym, amsfonts}
\def\@ptsize{2}
\newtheorem{Thm}{Theorem}

\newtheorem{Lem}{Lemma}

\newtheorem{Cor}{Corollary}

\newtheorem{Prop}{Proposition}
\newcommand{\qed}{\hfill $\Box$}

\begin{document}
\date{}
\maketitle {\def\thefootnote{} 
 {\def\thefootnote{}
}
\begin{abstract}
\noindent Robin's criterion states that the Riemann Hypothesis (RH) is true if and only
if Robin's inequality $\sigma(n):=\sum_{d|n}d<e^{\gamma}n\log \log n$ 
is satisfied for $n\ge  5041$, where 
$\gamma$ denotes the Euler(-Mascheroni) constant.
We show by elementary methods that if $n\ge 37$ does not satisfy Robin's criterion it
must be even and is neither squarefree nor
squarefull. Using a bound of Rosser and Schoenfeld we
show, moreover, that $n$ must be divisible by a fifth power $>1$. As 
consequence we obtain that RH holds true iff every natural number 
divisible by a fifth power $>1$ satisfies Robin's inequality. 
\end{abstract}
\section{Introduction}
Let $\cal R$ be the set of integers $n\ge 1$ satisfying $\sigma(n)<e^{\gamma}n\log \log n$. This 
inequality
we will call {\it Robin's inequality}. Note that it can be rewritten as
$$\sum_{d|n}{1\over d}<e^{\gamma}\log \log n.$$
 
Ramanujan \cite{Ra3} (in his original version
of his paper on highly composite integers, only part
of which, due to paper shortage, was published, for
the shortened version see \cite[pp. 78-128]{Ra}) proved that if RH holds then every sufficiently large
integer is in $\cal R$.
Robin \cite{R} proved that if RH holds, then actually every 
integer $n\ge 5041$
is in $\cal R$. He also showed that if RH is false, then there are infinitely many 
integers that are
not in $\cal R$. Put 
${\cal A}=\{1,2, 3, 4, 5, 6, 8, 9, 10, 12, 16, 18, 20, 24, 30, 36, 48, 60, 72, 84, 120, 180, 240,\hfil\break 360, 720, 
840, 2520,5040\}$. 
The set $\cal A$ consists of the integers $n\le 5040$ that do not satisfy Robin's inequality.
Note that none of the integers in $\cal A$ is divisible by a $5$th power of a prime.\\
\indent In this paper we are interested in establishing the inclusion of various infinite
subsets of the natural numbers in $\cal R$. We will
prove in this direction:
\begin{Thm}
\label{sqf}
Put ${\cal B}=\{2,3,5,6,10,30\}$. Every squarefree integer that is not in $\cal B$ is 
an element of $\cal R$.
\end{Thm} 
A similar result for the odd integers will be established:
\begin{Thm}
\label{odd}
Any odd positive integer $n$ distinct from $1,3,5$ and $9$ is in $\cal R$.
\end{Thm}
On combining Robin's result with the above theorems one finds:
\begin{Thm} The RH is true if and only for all even non-squarefree integers $\ge 5044$ 
Robin's inequality is satisfied.
\end{Thm}
It is an easy exercise to show that the even non-squarefree integers have density
${1\over 2}-{2\over \pi^2}=0.2973\cdots$ (cf. Tenenbaum \cite[p. 46]{T}).
Thus, to wit, this paper gives at least half a proof of RH !\\
\indent Somewhat remarkably perhaps these two results will be proved using only very
elementary methods. The deepest input will be Lemma \ref{beginner} below
which only requires pre-Prime Number Theorem elementary methods for its
proof (in Tenenbaum's \cite{T} introductory book on analytic number theory it is
already derived within the first 18 pages).\\
\indent Using a bound of Rosser and Schoenfeld (Lemma \ref{expmert} below), which ultimately
relies on some explicit knowledge regarding the first so many zeros of the 
Riemann zeta-function, one can prove some further results:
\begin{Thm}
\label{squarie}
The only squarefull integers not in $\cal R$ are $1,~4,~8,~9,~16$ and $36$. 
\end{Thm}
We recall that an integer $n$ is said to be squarefull if for every prime
divisor $p$ of $n$ we have $p^2|n$. An integer $n$ is called $t$-free if
$p^t\nmid m$ for every prime number $p$. 
(Thus saying a number is squarefree is the same as saying that
it is $2$-free.)
\begin{Thm}
\label{fivefree}
All {\rm 5}-free integers satisfy Robin's inequality.
\end{Thm}
Together with the observation that all exceptions $\le 5040$ to Robin's inequality
are 5-free and Robin's criterion, this result implies the following alternative 
variant of Robin's criterion.
\begin{Thm}
The RH holds iff for all integers $n$ divisible by the fifth power of some prime we have
$\sigma(n)<e^{\gamma}n\log \log n$.
\end{Thm}
The latter result has the charm of not involving a finite range of integers that
has to be excluded (the range $n\le 5040$ in Robin's criterion). We note that
a result in this spirit has been earlier established by 
Lagarias \cite{L} who, using
Robin's work, showed that the RH is equivalent with the inequality
$$\sigma(n)\le h(n)+e^{h(n)}\log(h(n)),$$
where $h(n)=\sum_{k=1}^n 1/k$ is the harmonic sum. 
\vfil\eject

\section{Proof of Theorem 1 and Theorem 2}
Our proof of Theorem \ref{sqf} requires the following lemmata.
\begin{Lem}$~$\\
\label{beginner}
{\rm 1)} For $x\ge 2$ we have 
$$\sum_{p\le x}{1\over p}=\log \log x + B + O({1\over \log x}),$$
where the implicit constant in Landau's symbol does not exceed
$2(1+\log 4)<5$ and 
$$B=\gamma + \sum_p \left(\log(1-{1\over p})+{1\over p}\right)=0.2614972128\cdots$$
denotes the (Meissel-)Mertens constant.\\
{\rm 2)} For $x\ge 5$ we have
$$\sum_{p\le x}{1\over p}\le \log \log x + \gamma.$$
\end{Lem}
{\it Proof}. 1) This result can be proved with very elementary methods. 
It is derived from scratch in the book of Tenenbaum \cite{T}, p. 16. At p. 18 the
constant $B$ is determined.\\
2) One checks that the inequality holds true for all primes $p$ satisfying
$5\le p\le 3673337$. On noting that
$$B+{2(1+\log 4)\over \log 3673337}<\gamma,$$
the result then follows from part 1. \qed\\

\noindent {\tt Remark 1}. More information on the (Meissel-)Mertens constant can be found e.g. in
the book of Finch \cite[\S 2.2]{F}.\\
{\tt Remark 2}. Using deeper methods from (computational) prime number theory Lemma \ref{beginner} can
be considerably sharpened, see e.g. \cite{RS}, but the point we want to make here is that
the estimate given in part 2, which is the estimate we need in the sequel, is a rather elementary
estimate.\\

We point out that $15$ is in $\cal R$.
\begin{Lem}
\label{prelude}
If $r$ is in $\cal A$ and $q\ge 7$ is a prime, then $rq$ is in $\cal R$, except when
$q=7$ and $r=12,~120$ or $360$.
\end{Lem}
\begin{Cor}
If $r$ is in $\cal B$ and $q\ge 7$ is a prime, then $rq$ is in $\cal R$.
\end{Cor}
{\it Proof of Lemma} \ref{prelude}. One verifies the result
in case $q=7$. Suppose that $r$ is in $\cal A$. 
Direct computation shows that $11r$ is in $\cal R$. {}From this
we obtain for $q\ge 11$ that
$${\sigma(rq)\over rq}=(1+{1\over q}){\sigma(r)\over r} \le {12\sigma(r)\over 11r}
< e^{\gamma}\log \log (11r)\le e^{\gamma}\log \log (qr).$$

\noindent {\it Proof of Theorem }\ref{sqf}. By induction with respect to $\omega(n)$, that
is the number of distinct prime factors of $n$. 
Put $\omega(n)=m$.
The assertion is easily provable for those integers with $m=1$ (the 
primes that is). Suppose it is
true for $m-1$, with $m\ge 2$ and let us consider the assertion for those squarefree $n$ with
$\omega(n)=m$. So let $n=q_1\cdots q_m$ be a squarefree number that is not in $\cal B$ and
assume w.l.o.g. that $q_1<\cdots < q_m$.
We consider two cases:\\
\noindent {\bf Case 1}: $q_m\ge \log (q_1\cdots q_m)=\log n$.\\ 
If $q_1\cdots q_{m-1}$ is in $\cal B$, then if $q_m$ is not in $\cal B$, $n=q_1\ldots
q_{m-1}q_m$ is in $\cal R$ (by the corollary to Lemma \ref{prelude}) and we are done, and if $q_m$ is
in $\cal B$, the only possibility is $n=15$ which is in $\cal R$ and we are
also done.
\newline If $q_1\cdots q_{m-1}$ is not in $\cal B$, by the induction
hypothesis we have 
$$(q_1+1)\cdots (q_{m-1}+1)< e^{\gamma}q_1\cdots q_{m-1}\log\log (q_1\cdots q_{m-1}),$$
and hence
\begin{equation}
\label{nogzeg}
(q_1+1)\cdots (q_{m-1}+1)(q_m+1)< e^{\gamma}q_1\cdots q_{m-1}(q_m+1)\log\log (q_1\cdots q_{m-1}).
\end{equation} 
We want to show that
$$
e^{\gamma}q_1\cdots q_{m-1}(q_m+1)\log\log (q_1\cdots q_{m-1})$$
\begin{equation}
\label{nogzeg1}
\leq e^{\gamma}q_1\cdots q_{m-1}q_m\log\log (q_1\cdots q_{m-1}q_m)=e^{\gamma}n\log \log n.
\end{equation}
Indeed (\ref{nogzeg1}) is equivalent with
$q_m\log\log (q_1\cdots q_{m-1}q_m)\geq (q_m+1)\log\log (q_1\cdots
q_{m-1}),$ or alternatively
\begin{equation}
\label{nogeen}
{q_m(\log\log (q_1\cdots
q_{m-1}q_m )-\log\log (q_1\cdots q_{m-1}))\over \log q_m}\ge 
{\log\log (q_1\cdots q_{m-1})\over \log q_m}.
\end{equation}
Suppose that $0<a<b$. Note that we have
\begin{equation}
\label{integralie}
{\log b - \log a\over b-a}= {1\over b-a}\int_a^b{dt\over t}> {1\over b}.
\end{equation}
Using this inequality we infer that (\ref{nogeen}) (and thus (\ref{nogzeg1})) is certainly satisfied 
if the next inequality
is satisfied:
$${q_m\over \log(q_1\cdots q_m)}\geq {\log\log(q_1\cdots
q_{m-1})\over \log q_m}.$$
Note that our assumption that $q_m\ge \log (q_1\cdots q_m)$ implies
that the latter inequality is indeed satisfied.\\
\noindent {\bf Case 2}: $q_m<\log (q_1\cdots q_m)=\log n$.\\ 
It is easy to see that $\sigma(n) < e^{\gamma}n\log \log n$ is equivalent
with
\begin{equation}
\label{lichie}
\log(q_1+1)-\log q_1 + \cdots + \log(q_m+1)-\log q_m < \gamma + \log \log \log 
(q_1\cdots q_m).
\end{equation}
Note that
$$\log(q_1+1)-\log q_1=\int_{q_1}^{q_1+1}{dt\over t}< {1\over q_1}.$$
In order to prove (\ref{lichie}) it is thus enough to prove that
\begin{equation}
\label{enough}
{1\over q_1}+\cdots + {1\over q_m}\le \sum_{p\le q_m}
{1\over p}\le \gamma + \log \log \log 
(q_1\cdots q_m).
\end{equation}
Since $q_m\ge 7$ we have by part 2 of Lemma \ref{beginner} and the 
assumption $q_m<\log (q_1\cdots q_m)$ that 
$$\sum_{p\le q_m}
{1\over p}\le \gamma + \log \log q_m < \gamma +  \log \log \log 
(q_1\cdots q_m),$$
and hence (\ref{enough}) is indeed satisfied. \qed\\

\noindent Theorem \ref{odd} will be derived from the following stronger result.
\begin{Thm}
\label{oddie}
For all odd integers except $1,3,5,9$ and $15$ we have
\begin{equation}
\label{nikkie2} 
{n\over \varphi(n)}<e^{\gamma}\log \log n,
\end{equation}
where $\varphi(n)$ denotes Euler's totient function.
\end{Thm}
To see that this is a stronger result, let $n=\prod_{i=1}^k p_i^{e_i}$ be the 
prime factorisation of $n$ and note that for $n\ge 2$ we have
\begin{equation}
\label{opstart}
{\sigma(n)\over n}=\prod_{i=1}^k{1-p_i^{-e_i-1}\over 1-p_i^{-1}}
<\prod_{i=1}^k{1\over 1-p_i^{-1}}={n\over \varphi(n)}.
\end{equation}

\noindent We let $\cal N$ ($\cal N$ in acknowledgement of
the contributions of J.-L. Nicolas to this subject) denote the set of integers $n\ge 1$ 
satisfying (\ref{nikkie2}). 
\noindent 
Our proofs of Theorems \ref{odd} and
\ref{oddie} use the next lemma.
\begin{Lem}
\label{specialset}
Put 
${\cal S}=\{3^a\cdot 5^b\cdot q^c:q\ge 7{\rm ~is~prime},~a,b,c\ge 0\}$. All elements from
$\cal S$ except $1,~3,~5$ and $9$ are in $\cal R$. All elements from $\cal S$ except
$1,~3,~5,~9$ and $15$ are in $\cal N$.
\end{Lem}
{\it Proof}. If $n$ is in $\cal S$ and $n\ge 31$ we have
$${\sigma(n)\over n}\le {n\over \varphi(n)}\le 
{3\over 2}\cdot {5\over 4}\cdot {q\over q-1}\le 
{3\over 2}\cdot {5\over 4}\cdot {7\over 6}<e^{\gamma} \log \log n.$$
Using this observation the proof is easily completed. \qed\\

\noindent {\tt Remark}. Let $y$ be any integer. Suppose that we have an infinite set of integers all having
no prime factors $>y$. Then $\sigma(n)/n$ and $n/\varphi(n)$ 
are bounded above on this set, whereas
$\log \log n$ tends to infinity. Thus only finitely many of those integers will not 
be in $\cal R$, respectively 
$\cal N$.
It is a finite computation to find them all (cf. the proof of Lemma \ref{specialset}).\\

\noindent {\it Proof of Theorem} \ref{oddie}. 
As before we let $m=\omega(n)$. If $m\le 1$ then, by Lemma \ref{specialset}, $n$ is in $\cal N$, except 
when $n=1,3,5$ or $9$. So we may assume $m\ge 2$.
Let $\kappa(n)=\prod_{p|n}p$ denote the squarefree kernel of
$n$. Since $n/\varphi(n)=\kappa(n)/\varphi(\kappa(n))$ it follows that if $r$ is a squarefree
 number satisfying (\ref{nikkie2}), then all integers $n$ with $\kappa(n)=r$ satisfy (\ref{nikkie2}) 
 as well. Thus we consider first the case where $n=q_1\cdots q_m$ is an odd  squarefree integer with
 $q_1< \cdots < q_m$.
In this case $n$ is in $\cal N$ iff
$${n\over \varphi(n)}=\prod_{i=1}^m{q_i\over q_i-1}< e^{\gamma}\log \log n.$$
Note that 
$${q_i\over q_i-1}\le {3\over 2}{\rm ~and~}{q_i\over q_i-1}< {q_{i-1}+1\over q_{i-1}},$$
and hence
$${n\over \varphi(n)}=\prod_{i=1}^m{q_i\over q_i-1}< {3\over 2}\prod_{i=1}^{m-1}{q_i+1\over q_i}
={\sigma(n_1)\over n_1},$$
where $n_1=2n/q_m<n$. Thus, $n/\varphi(n)<\sigma(n_1)/n_1$.
If $n_1$ is in $\cal R$, then invoking Theorem \ref{sqf} we find
$${n\over \varphi(n)}< {\sigma(n_1)\over n_1} < e^{\gamma}\log \log n_1 < 
e^{\gamma} \log \log n,$$
and we are done.\\
\indent If $n_1$ is not in $\cal R$, then by Theorem \ref{sqf} it follows that $n$ must be in $\cal S$.
The proof is now completed on invoking Lemma \ref{specialset}. \qed\\

\noindent {\it Proof of Theorem} \ref{odd}. 
One checks that $1,3,5$ and $9$ are not in $\cal R$, but $15$ is in $\cal R$.
The result now follows by Theorem \ref{oddie} and inequality (\ref{opstart}). \qed

\subsection{Theorem \ref{oddie} put into perspective}
Since the proof of Theorem \ref{oddie} can be carried out with such simple means, one might
expect it can be extended to quite a large class of even integers. However, even a superficial
inspection of the literature on $n/\varphi(n)$ shows this expectation to be wrong.\\
\indent Rosser and Schoenfeld \cite{RS} showed in 1962 that
$${n\over \varphi(n)}\le e^{\gamma}\log \log n + {5\over 2\log \log
n},$$ with one exception: $n=2\cdot 3\cdot 5\cdot 7\cdot 11\cdot 13\cdot 
17\cdot 19\cdot 23$. 
They raised the question of whether
there are infinitely many $n$ for which
\begin{equation}
\label{nikkie} 
{n\over \varphi(n)}>e^{\gamma}\log \log n, 
\end{equation}
which was answered in the
affirmative by J.-L. Nicolas \cite{N}. More precisely, let $N_k=2\cdot 3\cdot \cdots p_k$ be the product
of the first $k$ primes, then if the RH holds true (\ref{nikkie}) is satisfied 
with $n=N_k$ for every $k\ge 1$. On the other hand, if RH is false, then there
are infinitely many $k$ for which (\ref{nikkie}) is satisfied with $n=N_k$ and there 
are infinitely many $k$ for which (\ref{nikkie}) is not satisfied
with $n=N_k$. Thus the approach we have taken to prove Theorem \ref{odd}, namely to derive it
from the stronger result Theorem \ref{oddie}, is not going to work for even integers.

\section{Proof of Theorem \ref{squarie}}
The proof of Theorem \ref{squarie} is an immediate consequence
of the following stronger result.
\begin{Thm}
\label{squarie2}
The only squarefull integers $n\geq 2$ not in
$\cal N$ are $4,~ 8,~ 9,~ 16,~ 36,~ 72,~ 108$, $144,~ 216,~
900,~ 1800,~
2700,~ 3600,~ 44100$ and $88200.$
\end{Thm}
Its proof requires the following two lemmas.
\begin{Lem} {\rm \cite{RS}}.
\label{expmert}
For $x>1$ we have
$$\prod_{p\le x}{p\over p-1}\le e^{\gamma}(\log x+{1\over \log x}).$$
\end{Lem}
\begin{Lem}
\label{ddrie}
Let $p_1=2,~p_2=3,\ldots$ denote the consecutive primes. If 
$$\prod_{i=1}^m {p_i \over {p_i -1}} \ge e^{\gamma} \log (2\log (p_1\cdots p_m)),$$
then $m\leq 4.$
\end{Lem}
{\it Proof}. Suppose that $m\ge 26$ (i.e. $p_m\ge 101$). It then
follows by Theorem 10 of \cite{RS}, which 
states that $\theta(x):=\sum_{p\le x}\log p>0.84x$
for $x\ge 101$, that $\log(p_1\cdots p_m)=\theta(p_m)>0.84p_m$. We find that
$$\log(2\log(p_1\cdots p_m))>\log p_m+\log 1.64\ge \log p_m +{1\over \log p_m},$$
and so, by Lemma \ref{expmert}, that 
$$\prod_{i=1}^m{p_i\over p_i-1}\le e^{\gamma}\left(\log p_m+{1\over \log p_m}\right)
<e^{\gamma}\log(2\log(p_1\cdots p_m)).$$
The proof is then completed on checking the inequality directly for the
remaining values of $m$. \qed\\

\noindent {\it Proof of Theorem} \ref{squarie2}. Suppose that 
$${n\over {{\varphi}(n)}} \ge e^{\gamma} \log \log n.$$
Put $\omega (n)=m.$ Then
$$\prod_{i=1}^m {{p_i}\over {p_i -1}} \geq {n\over {\varphi (n)}} \geq e^{\gamma}
\log \log n \geq e^{\gamma}
\log (2 \log (p_1~\cdots~p_n)).$$
By Lemma \ref{ddrie} it follows that $m\leq 4.$ In particular we must have
$$2\cdot {3\over 2}\cdot {5\over 4}\cdot {7\over 6} = {35\over 8} \geq e^{\gamma}
\log \log n,$$
whence $n\leq \exp (\exp ( e^{-\gamma} 35/8)) \leq 116144.$
On numerically checking the inequality for the squarefull integers $\le 116144$,
the proof is then completed. \qed\\

\noindent {\tt Remark}. The squarefull integers $\le 116144$ are easily 
produced on noting
that they can be unqiuely written as $a^2b^3$, with $a$ a positive integer and $b$ squarefree.

\section{On the ratio $\sigma(n)/(n\log \log n)$ as $n$ ranges over various sets
of integers}
We have proved that Robin's inequality holds for large enough odd numbers, squarefree 
and squarefull numbers.
A natural question to ask is how large the ratio $f_1(n):=\sigma(n)/(n\log \log n)$ can be when
we restrict $n$ to these sets of integers. We will consider the same question for the
ratio $f_2(n):=n/(\varphi(n)\log \log n)$. Our results in this direction are summarized
in the following result:
\begin{Thm}
\label{6}
We have
$${\rm (1)~}\limsup_{n\rightarrow \infty}f_1(n)={e^{\gamma}},~
{\rm (2)~}\limsup_{n\rightarrow \infty\atop n{\rm ~is~squarefree}}f_1(n)={6e^{\gamma}\over \pi^2},~
{\rm (3)~}\limsup_{n\rightarrow \infty\atop n {\rm ~is~odd}}f_1(n)={e^{\gamma}\over 2},$$
and, moreover,
$${\rm (4)~}\limsup_{n\rightarrow \infty}f_2(n)={e^{\gamma}},~
{\rm (5)~}\limsup_{n\rightarrow \infty\atop n{\rm ~is~squarefree}}f_2(n)={e^{\gamma}},~
{\rm (6)~}\limsup_{n\rightarrow \infty\atop n {\rm ~is~odd}}f_2(n)={e^{\gamma}\over 2}.$$
Furthermore,
$${\rm (7)~}\limsup_{n\rightarrow \infty\atop n{\rm ~is~squarefull}}f_1(n)={e^{\gamma}},~
{\rm (8)~}\limsup_{n\rightarrow \infty\atop n{\rm ~is~squarefull}}f_2(n)=e^{\gamma}.$$
\end{Thm}
(The fact that the corresponding lim infs are all zero is immediate on letting $n$ run
over the primes.)\\
\indent Part 4 of Theorem \ref{6} was proved by Landau in 1909, see e.g. \cite[Theorem 13.14]{A}, and
the remaining parts can be proved in a similar way. Gronwall in 1913 established part 1.
Our proof makes use of a lemma involving
$t$-free integers (Lemma \ref{pro}), which is easily proved on invoking 
a celebrated result due to Mertens (1874) asserting that
\begin{equation}
\label{franzm}
\prod_{p\le x}\left(1-{1\over p}\right)^{-1}\sim e^{\gamma} \log x,~~x\rightarrow \infty.
\end{equation}
\begin{Lem}
\label{pro}
Let $t\ge 2$ be a fixed integer.
We have
$${\rm (1)~}\limsup_{n\rightarrow \infty\atop t{\rm -free~integers}}f_1(n)={e^{\gamma}\over \zeta(t)},~
{\rm (2)~}\limsup_{n\rightarrow \infty\atop {\rm odd~}t{\rm -free~integers}}f_1(n)={e^{\gamma}\over 2\zeta(t)(1-2^{-t})}.$$
\end{Lem}
{\it Proof}. 1) Let us consider separately the prime divisors of $n$ that are larger
than $\log n$. Let us say there are $r$ of them. Then $(\log n)^r<n$ and thus
$r<\log n/\log \log n$. Moreover, for $p>\log n$ we have
$${1-p^{-t}\over 1-p^{-1}}<{1-(\log n)^{-t}\over 1-(\log n)^{-1}}.$$
Thus,
$$\prod_{p|n\atop p>\log n}{1-p^{-t}\over 1-p^{-1}}
< \left({1-(\log n)^{-t}\over 1-(\log n)^{-1}}\right)^{\log n\over
\log \log n}.$$
Let $p_k$ denote the largest prime factor of $n$. We obtain
\begin{eqnarray}
\label{ketting}
{\sigma(n)\over n}&=&\prod_{i=1}^k{1-p_i^{-e_i-1}\over 1-p_i^{-1}}
\le \prod_{i=1}^k{1-p_i^{-t}\over 1-p_i^{-1}}\cr
&<&
\left({1-(\log n)^{-t}\over 1-(\log n)^{-1}}\right)^{\log n\over
\log \log n}\prod_{p\le \log n}{1-p^{-t}\over 1-p^{-1}},
\end{eqnarray}
where in the derivation of the first inequality we used that $e_i<t$ by
assumption. Note that the factor before the final product
satisfies $1+O((\log \log n)^{-1})$ and thus tends to $1$ as $n$ tends
to infinity.
On invoking (\ref{franzm}) and 
noting that $\prod_{p\le \log n}(1-p^{-t})\sim \zeta(t)^{-1}$, it
follows that the $\limsup \le e^{\gamma}/\zeta(t)$.\\
\indent In order to prove the $\ge$ part of the assertion, take $n=\prod_{p\le x}p^{t-1}$. Note
that $n$ is $t$-free.
On invoking (\ref{franzm}) 
we infer that 
$${\sigma(n)\over n}=\prod_{p\le x}{1-p^{-t}\over 1-p^{-1}}\sim {e^{\gamma}\over \zeta(t)}
\log x.$$
Note that $\log n=t\sum_{p\le x}p=t\theta(x)$, 
where $\theta(x)$ denotes the Chebyshev theta function. By an equivalent form
of the Prime Number Theorem we have $\theta(x)\sim x$ and hence
$\log \log n=(1+o_t(1))\log x$. It follows that for the particular sequence
of infinitely many $n$ values under consideration we have
$${\sigma(n)\over n\log \log n}={e^{\gamma}\over \zeta(t)}\biggl(1+o_t(1)\biggr).$$ 
Thus, in particular, for a given $\epsilon>0$ there are 
infinitely many $n$  such that 
$${\sigma(n)\over n\log \log n}>{e^{\gamma}\over \zeta(t)}(1-\epsilon).$$
2) Can be proved very similarly to part 1. Namely, the third product in (\ref{ketting}) 
will extend over the primes $2< p\le \log n$ and for the $\ge $ part we consider the
integers $n$ of the form $n=\prod_{2<p\le x}p^{t-1}$. 
\qed\\

\noindent {\tt Remark}. Robin \cite{R} has shown that if RH is false, then
there are infinitely many integers $n$ not in $\cal R$. As 
$n$ ranges over these numbers, then by part 1 of Lemma \ref{pro} we must
have $\max\{e_i\}\rightarrow \infty$, 
where $n=\prod_{i=1}^k p_i^{e_i}$.\\

\noindent {\it Proof of Theorem} \ref{6}.\\
{\bf 1)} Follows from part 1 of Lemma \ref{pro} on letting $t$
tend to infinity. A direct proof (similar to that of Lemma \ref{pro}) can also be
given, see e.g. \cite{B}. This result was proved first by Gronwall in 1913.\\
{\bf 2)} Follows from part 1 of Lemma \ref{pro} with $t=2$.\\
{\bf 3)} Follows on letting $t$ tend to infinity in part 2 of Lemma \ref{pro}.\\
{\bf 4)} Landau (1909).\\
{\bf 5)} Since $f_2(n)\le f_2(\kappa(n))$, part 5 is a consequence of part 4.\\
{\bf 6)} A consequence of part 4 and the fact that for odd integers $n$ and $a\ge 1$ we have
$f_2(2^a n)=2f_2(n)(1+O((\log n\log \log n)^{-1}))$.\\
{\bf 7)} Consider numbers of the form $n=\prod_{p\le x}p^{t-1}$ and let $t$ tend to
infinity. These are squarefull for $t\ge 3$ and using them the $\ge$ part of the assertion
follows. The $\le $ part follows of course from part 3.\\
{\bf 8)} It is enough here to consider the squarefull numbers of the form $n=\prod_{p\le x}p^2$. \qed

\section{Reduction to Hardy-Ramanujan integers}
Recall that $p_1,p_2,\ldots$ denote the consecutive primes. An integer of
the form $\prod_{i=1}^s p_i^{e_i}$ with $e_1\ge e_2\ge \cdots \ge e_s\ge 0$ we will
call an {\it Hardy-Ramanujan integer}. We name them after Hardy and Ramanujan who
in a paper entitled `A problem in the analytic theory of numbers' (Proc. London Math.
Soc. {\bf 16} (1917), 112-132) investigated them. See also \cite[pp. 241-261]{Ra}, where
this paper is retitled `Asymptotic formulae for the distribution of integers of various
types'.
\begin{Prop}
\label{HRnumber}
If Robin's inequality holds for all Hardy-Ramanujan integers $5041\le n\le x$, then
it holds for all integers $5041\le n\le x$.
Asymptotically there are $$\exp((1+o(1))2\pi\sqrt{\log x/3\log \log x})$$ Hardy-Ramanujan numbers
$\le x$.
\end{Prop}
Hardy and Ramanujan proved the asymptotic assertion above.
The proof of the first part requires a few lemmas.
\begin{Lem}
\label{flauwie}
For $e>f>0,$ the function 
$$g_{e,f}:x\to {{1-x^{-e}}\over {1-x^{-f}}}$$ is
strictly decreasing on $(1,+\infty ]$.
\end{Lem}
{\it Proof}. For $x>1,$ we have $$g_{e,f}'(x)={{ex^f-fx^e+f-e}\over
{x^{e+f+1}{\left( 1-x^{-f}\right) }^2}}.$$ Let us consider the function
$h_{e,f}:x\to ex^f-fx^e+f-e$. For $x>1,$ we have 
$h_{e,f}'(x)=efx^f\left( 1-x^{e-f}\right) <0$.
Consequently $h_{e,f}$ is decreasing on $(1,+\infty]$ and since $h_{e,f}(1)=0,$
we deduce that $h_{e,f}(x)<0$ for $x>1$ and so $g_{e,f}(x)$ is strictly
decreasing on $(1,+\infty]$. \qed\\

\noindent {\tt Remark}. In case $f$ divides $e,$ then 
$${{1-x^{-e}}\over {1-x^{-f}}}= 1+ {1\over {x^f}}+{1\over {x^{2f}}}+\cdots +{1\over
{x^e}},$$
and the result is obvious.
\begin{Lem}
\label{eenhoorn}
If $q>p$ are primes and $f>e,$ then
\begin{equation}
\label{star}
{{\sigma\left( p^fq^e\right)}\over {p^fq^e}} > {{\sigma \left(p^eq^f\right) }\over
{p^eq^f}}. 
\end{equation}
\end{Lem}
{\it Proof.} Note that the inequality (\ref{star}) is equivalent with
$$(1-p^{-1-f})(1-p^{-1-e})^{-1} > 
(1-q^{-1-f})(1-q^{-1-e})^{-1}.$$
It follows by Lemma \ref{flauwie} that the latter inequality is satisfied.
\qed\\

Let $n=\prod_{i=1}^s {q_i}^{e_i}$ be a factorisation of $n,$ where we ordered
the primes $q_i$ in such a way that $e_1 \geq e_2 \geq e_3\ge \cdots$ We say that
${\bar e}=(e_1,\ldots,e_s)$ is the exponent pattern of the integer $n$. Note that
$\Omega (n)=e_1+\ldots+e_s$, where $\Omega (n)$ denotes the total number of
prime divisors of $n$. Note that $\prod_{i=1}^s {p_i}^{e_i}$ is the minimal
number having exponent pattern ${\bar e}$. We denote this 
(Hardy-Ramanujan) number by $m({\bar e})$.

\begin{Lem}
\label{vvier}
We have
$$\max \left\{ {{\sigma(n)}\over n} ~|~ n \mbox{ has factorisation pattern } {\bar e}
\right\} =
{\sigma(m({\bar e}))\over m({\bar e})}.$$
\end{Lem}
{\it Proof}. Since clearly $\sigma(p^e)/p^e > \sigma(q^e)/q^e$
if $p<q,$ the maximum is assumed on integers $n=\prod_{i=1}^s {p_i}^{f_i}$
having factorisation pattern ${\bar e}.$ Suppose that $n$ is any number of this
form for which the maximum is assumed, then by Lemma \ref{eenhoorn} it follows that $f_1 \geq
f_2 \geq \cdots \geq f_s$ and so $n=m({\bar e})$. \qed

\begin{Lem}
\label{vvijf}
Let ${\bar e}$ denote the factorisation pattern of $n$.\\
{\rm 1)} If ${{\sigma(n)}/n} \geq e^{\gamma} \log \log n,$ then ${{\sigma(m({\bar
e}))}/{m({\bar e})}} \ge e^{\gamma} \log \log m({\bar e}).$\\
{\rm 2)} If ${{\sigma(m({\bar e}))}/{m({\bar e})}} < e^{\gamma} \log \log
m({\bar e}),$ then ${{\sigma(n)}/n} < e^{\gamma} \log \log n$ for every integer
$n$ having exponent pattern ${\bar e}$.
\end{Lem}
{\it Proof}. A direct consequence of the fact that $m({\bar e})$ is the smallest
number having exponent pattern ${\bar e}$ and Lemma \ref{vvier}. \qed\\

On invoking the second part of the latter lemma, the proof of Proposition \ref{HRnumber}
is completed.

\section{The proof of Theorem \ref{fivefree}}
Our proof of Theorem \ref{fivefree} makes use of lemmas \ref{gladjakker}, \ref{doenbaar}
and \ref{proddie}.
\begin{Lem}
\label{gladjakker}
Let $t\ge 2$ be fixed.
Suppose that there exists a $t$-free integer exceeding $5040$ that does
not satisfy Robin's inequality. Let $n$ be the smallest 
such integer. Then $P(n)<\log n$, where $P(n)$ denotes the largest prime factor of $n$.
\end{Lem}
{\it Proof}. Write $n=r\cdot q$ with $P(n)=q$ and note that $r$ is $t$-free. The minimality
assumption on $n$ implies that either $r\le 5040$ and does not satisfy Robin's inequality or 
that $r$ is in $\cal R$.
First assume we are in the former case. Since 720 is the largest integer $a$ in $\cal A$ with
$P(a)\le 5$ and $5\cdot 720\le 5040$, it follows that $q\ge 7$. By Lemma \ref{specialset} we then 
infer, using the assumption that $n>5040$,  that $n=qr$ is in $\cal R$; a contradiction.
Thus we
may assume that $r$ is in $\cal R$ and therefore $r\ge 7$. We will now show that this together with the
assumption $q\ge \log n$ leads to a contradiction, whence the result
follows.

So assume that $q\ge \log n$. This implies that $q\log q\ge \log n\log \log n>\log n\log \log r$
and hence
$${q\over \log n}>{\log \log r\over \log q}.$$
This  implies that
\begin{equation}
\label{burp}
{q(\log\log n-\log \log r)\over \log q}>{\log \log r\over \log q},
\end{equation}
where we used that
$${\log \log n-\log \log r\over \log q}={1\over \log n-\log r}\int_{\log r}^{\log n}{dt\over t}>{1\over \log n}.$$
Inequality (\ref{burp}) is equivalent with
$(1+{1/q})\log \log r < \log \log n$. Now we infer that 
\begin{equation}
\label{geblaat}
{\sigma(n)\over n}={\sigma(qr)\over qr}\le \left(1+{1\over q}\right){\sigma(r)\over r}
<\left(1+{1\over q}\right)e^{\gamma}\log \log r< e^{\gamma}\log \log n,
\end{equation}
where we used that $\sigma$ is submultiplicative (that is  $\sigma(qr)\le \sigma(q)\sigma(r)$).
The inequality (\ref{geblaat}) contradicts our assumption that $n\not\in \cal R$. \qed

\begin{Lem}
\label{doenbaar}
All 5-free Hardy-Ramanujan integers $n>5040$ with $P(n)\le 73$ satisfy Robin's inequality.
\end{Lem}
{\it Proof}. There are 12649 $5$-free Hardy-Ramanujan integers $n$ with $P(n)\le 73$, that
are easily produced using MAPLE. A further MAPLE computation learns that all integers
exceeding 5040 amongst these (12614 in total) are in $\cal R$. \qed\\

\noindent {\tt Reamrk}. On noting that $\prod_{p\le 73}p^4<\prod_{p\le 20000}p$ and invoking
Robin's result \cite[p. 204]{R} that an integer $n\not\in \cal R$ with $n>5040$ satisfies
$n\ge \prod_{p\le 20000}p$, an alternative proof of Lemma \ref{doenbaar} is obtained.

\begin{Lem}
\label{proddie}
For $x\ge 3$ and $t\ge 2$ we have that 
$$\sum_{p\le x}\log\left({1-p^{-t}\over 1-p^{-1}}\right)\le -\log \zeta(t)+{t\over (t-1)}x^{1-t}+\gamma+\log \log x
+\log\left(1+{1\over \log ^2 x}\right).$$
\end{Lem}
The proof of this lemma on its turn rests on the lemma below.
\begin{Lem}
\label{rtl}
Put $R_t(x)=\prod_{p>x}(1-p^{-t})^{-1}$.
For $x\ge 3$ and $t\ge 2$ we have that $\log(R_t(x))\le {tx^{1-t}/(t-1)}$.
\end{Lem}
{\it Proof}.
We have
\begin{eqnarray}
R_t(x)&=& -\sum_{p>x}\log\left(1-{1\over p^t}\right)=\sum_{p>x}\sum_{m=1}^{\infty}{1\over mp^{tm}}
\le \sum_{p>x}\sum_{m=1}^{\infty}{1\over (p^m)^t}\nonumber\\
&\le & \sum_{n>x}{1\over n^t}\le {1\over x^t}+\sum_{n>x+1}{1\over n^t}
\le {1\over x^{t-1}}+\int_x^{\infty}{du\over u^t}={t\over t-1}x^{1-t}.\nonumber
\end{eqnarray}

\noindent {\it Proof of Lemma} \ref{proddie}. On noting that $\prod_{p\le x}(1-p^{-t})=R_t(x)/\zeta(t)$ and invoking 
Lemma \ref{rtl} we obtain 
$$\sum_{p\le x}\log\left(1-{1\over p^{t}}\right)=-\log \zeta(t) +\log (R_t(x))\le -\log \zeta(t) + 
{t\over t-1}x^{1-t}.$$
On combining this estimate with Lemma \ref{expmert}, the estimate then follows. \qed

\begin{Lem}
\label{puzzelstukje}
Let $m$ be a 5-free integer such that $P(m)<\log m$ and
$m$ does not satisfy Robin's inequality. Then $P(m)\le 73$.
\end{Lem}
{\it Proof}. Put $t=5$.
Write $P_t(x)=\prod_{p\le x}(1-p^{-t})/(1-p^{-1})$.
Put $z=\log m$. 
The assumptions on $m$ imply that $\sigma(m)/m\le P_t(z)$.
This inequality in combination with Lemma \ref{proddie} yields 
$$\log \left({\sigma(m)\over m}\right)\le -\log \zeta(t) +{t\over (t-1)z^{t-1}}+\gamma+\log \log z 
+\log\left(1+{1\over \log ^2 z}\right).$$
Once $$ -\log \zeta(t) +{t\over t-1}z^{1-t}+\gamma+\log \log z 
+\log\left(1+{1\over \log ^2 z}\right) < \gamma + \log \log z,$$ Robin's inequality is satisfied. 
We infer that once we have found a $z_0\ge 3$ such
that
$${t\over t-1}z_0^{1-t}+\log\left(1+{1\over \log^2 z_0}\right)-\log \zeta(t)<0,$$
then Robin's inequality will be satisfied in case $z\ge z_0$.
One finds that $z_0=196$ will do. It follows that
$z<196$ and hence ${\sigma(m)/m}<P_5(193)=9.18883221\ldots$. Note that if
$e^{\gamma}\log \log m\ge P_5(193)$, then Robin's inequality is satisfied.
We thus conclude that $\log m\le 
\exp(P_5(193)e^{-\gamma})=174.017694\ldots$. Since 173 is the largest prime $<175$ we
know that $m$ must satisfy ${\sigma(m)/m}<P_5(173)=8.992602079\ldots$. We now
proceed as before, but with $P_5(193)$ replaced by $P_5(173)$. Indeed, this `cascading down'
can be repeated several times before we cannot reduce further. This is at the point
where we have reached the conclusion that $z=\log m\le 73$. Then we cannot reduce
further since $\exp(P_5(73)e^{-\gamma})>73$. \qed\\

\noindent {\it Proof of Theorem} \ref{fivefree}. By contradiction. So suppose a 5-free integer
exceeding 5040 exists that does not satisfy Robin's inequality. We let $n$ be the smallest
such integer. By Lemma \ref{gladjakker} it follows that $P(n)<\log n$, whence
by Lemma \ref{puzzelstukje} we infer that $P(n)\le 73$. We will now show that
$n$ is a Hardy-Ramanujan number. On invoking Lemma \ref{doenbaar} the proof is then
completed.

It thus remains to establish that $n$ is a Hardy-Ramanujan number. Let $\bar e$ denote the
factorisation pattern of $n$. Note that $m({\bar e})$ is 5-free and that $m({\bar e})\le n$. By
the minimality of $n$ and part 1 of Lemma \ref{vvijf} it follows that we cannot have
that $5041\le m({\bar e}) < n$ and so either $m({\bar e})=n$, in which case we are done
as $m({\bar e})$ is a Hardy-Ramanujan number, or $m({\bar e})\le 5040$.
In the latter case
we must have $n=p_1^{e_1}p_2^{e_2}p_3^{e_3}p_4^{e_4}p_5^{e_5}$ (since max$\{\omega(r)~:~r\le 5040\}=5$)
and so
$${\sigma(n)\over n}\le \prod_{p\le 11}{1-p^{-5}\over 1-p^{-1}}=4.6411\cdots$$
and
$$\prod_{p\le 11}{1-p^{-5}\over 1-p^{-1}}\ge e^{\gamma}\log \log n,$$ whence $\log n\le 13.55$.
A MAPLE computation now shows that $n\in {\cal R}$, contradicting our assumption that $n\not \in \cal R$. \qed\\

\noindent By the method above we have not been able to replace 5-free
by 6-free in Theorem \ref{fivefree} (this turns out to require a substantial
computational effort). Recently J.-L. Nicolas kindly informed the authors of an approach (rather
different from the one followed here and being less self-contained) that might lead to a serious improvement
of the 5-free. It would certainly be interesting to pursue Nicolas's idea
further and this might be part of a follow-up paper.

\section{Acknowledgement}
We thank J.-C. Lagarias for pointing out reference \cite{KB}. Furthermore,
E. Bach, P. Dusart, O. Ramar\'e and, especially, J.-L. Nicolas, for their remarks.
Keith Briggs we thank for his willingness to do large scale computations
on our behalf. 
In the end, however, it turned out that only modest
computations are needed in order to establish Theorem 6 (the main result of
this paper).

\vfil\eject
\medskip\noindent {\footnotesize 
{\tt Y.-J.~Choie}, Dept of Mathematics, POSTECH, Pohang, Korea 790-784,~e-mail: {\tt yjc@postech.ac.kr}\\
{\tt P. Moree}, Max-Planck-Institut f\"ur Mathematik,
Vivatsgasse 7, D-53111 Bonn, Germany.\\ e-mail: {\tt moree@mpim-bonn.mpg.de}\\
{\tt N.~Lichiardopol}, ESSI, Route des Colles, 06 903 Sophia Antipolis, France,~e-mail: {\tt lichiard@essi.fr}\\  
{\tt P.~Sol\'e}, CNRS-I3S, ESSI,
Route des Colles,
06 903 Sophia Antipolis,
France,~e-mail: {\tt ps@essi.fr}}

\end{document}